\documentclass[a4paper,12pt]{amsart}
\usepackage{amssymb}

\begin{document}

\newcommand{\D}{\ensuremath{\mathcal{D}}}
\newcommand{\I}{\ensuremath{\mathcal{I}}}
\newcommand{\Z}{\mathbb{Z}}
\newcommand{\C}{\mathbb{C}}
\renewcommand{\P}{\mathbb{P}}
\newcommand{\SL}{\operatorname{SL}}
\renewcommand{\H}{\operatorname{H}}
\newcommand{\SHom}{\mathcal{H}om}
\newcommand{\F}{\ensuremath{\mathcal{F}}}
\renewcommand{\O}{\mathcal{O}}
\newcommand{\SH}{\mathcal{H}}
\newcommand{\M}{\mathcal{M}}
\newcommand{\LW}{\mathcal{L}(v)}
\renewcommand{\L}{\mathcal{L}}
\newcommand{\g}{\mathfrak{g}}
\newcommand{\height}{\operatorname{ht}}
\newcommand{\kar}{\operatorname{char}}
\newcommand{\Spec}{\operatorname{Spec}}
\newcommand{\Supp}{\operatorname{Supp}}
\newcommand{\Sing}{\operatorname{Sing}}
\newcommand{\Gr}{\operatorname{Gr}}
\newcommand{\cd}{\operatorname{cd}}
\newtheorem{Corollary}{Corollary}[section]
\newtheorem{Lemma}{Lemma}[section]
\newtheorem{Theorem}{Theorem}[section]
\newtheorem{Proposition}[Theorem]{Proposition}
\newtheorem{Definition}[Theorem]{Definition}
\newtheorem{Remark}{Remark}
\newenvironment{Proof}{{\it Proof.\/}}{\hfill $\square$\medskip}

\title[Global $F$-regularity and $\D$-modules]{Global $F$-regularity of Schubert varieties with applications to $\D$-modules}
\author[N.~Lauritzen, U.~Raben-Pedersen and J.~F.~Thomsen]{Niels Lauritzen, Ulf Raben-Pedersen \and \\Jesper Funch Thomsen}
\address{Institut for matematiske fag\\Aarhus Universitet\\ \AA rhus, Denmark}
\email{niels@imf.au.dk, ab061278@imf.au.dk, funch@imf.au.dk}

\maketitle

A projective algebraic variety $X$ over an algebraically closed field
$k$ of positive characteristic is called globally $F$-regular 
if the section ring $S(\L) = \oplus_{n\geq 0} \H^0(X, \L^n)$ of
an ample line bundle $\L$ on $X$ is strongly $F$-regular in
the sense of Hochster and Huneke \cite{HH} (cf.~Definition \ref{DefFreg}). 
This important
notion was introduced by Karen Smith in \cite{SmithFreg}.

In this paper we prove that Schubert varieties are globally$F$-regular. 
An immediate consequence is that local rings
of Schubert varieties are strongly $F$-regular and thereby $F$-rational.
Another consequence is that local rings of varieties (like determinantal
varieties) that can be identified with open subsets of
Schubert varieties (cf.~\cite{LM}) are strongly $F$-regular

Let $X$ denote a flag variety and $Y\subset X$ a Schubert variety over
$k$. 
Then the local cohomology sheaves $\SH^j_Y(\O_X)$ are equivariant
(for the action of the Borel subgroup) and holonomic (in the sense
of \cite{Bogvad02}) $\D_X$-modules.
As an application of $F$-rationality of Schubert varieties we apply 
recent results of 
Blickle (cf.~\cite{Blickle}) to
prove that the simple objects in the category of
equivariant and holonomic $\D_X$-modules are precisely
the local cohomology sheaves $\SH^c_Y(\O_X)$, where $c$ is
the codimension of $Y$ in $X$. Using a local  
Grothendieck-Cousin complex from \cite{KL}, we prove that the
decomposition of the local cohomology modules with support in
Bruhat cells is multiplicity free (see \S \ref{GC}). 

In characteristic zero the local cohomology modules with support in
Bruhat cells correspond to dual Verma modules. In this
setting the decomposition behavior and the simple $\D_X$-modules 
arise from intersection cohomology 
complexes of Schubert varieties by the Riemann-Hilbert correspondence.  
Picking the singular codimension one Schubert variety $Y$ in the full 
flag variety $Z$ for $\SL_4$ in characteristic zero,
computations in Kazhdan-Lusztig theory show that $\SH^1_Y(\O_Z)$
is not a simple $\D_Z$-module (see \S \ref{nonsmpl}). 

We are grateful to M.~Kashiwara and V.~B.~Mehta for discussions related
to this work. We also thank the referee for pointing out the
connection to the Riemann-Hilbert correspondence in positive
characteristic by Emerton and Kisin \cite{EK} in proving simplicity
of local cohomology.

\section{Global $F$-regularity and Frobenius splitting} 

Let $k$ denote a field of 
characteristic $p>0$ and $R$ a finitely 
generated $k$-algebra. For an $R$-module  $M$ we 
define $F^e_* M$ to be the $R$-module which is 
equal to $M$ as an abelian group with 
$R$-action given by $r \cdot m = r^{p^e} m$.  

\begin{Definition}(\cite{HH})\label{DefFreg}
The ring $R$ is said to be {\it strongly $F$-regular}
if for every $c \in R$, not contained in any minimal prime of 
$R$, there exists a positive integer $e \geq 0$ such 
that the map of $R$-modules 
$$R \rightarrow F_*^e R,$$
$$ 1 \mapsto c $$ 
is split.
\end {Definition}  

The concept of strong $F$-regularity has been extended 
to projective varieties over $k$ in the following way.

\begin{Definition}(\cite{SmithFreg})
The projective variety  $X$ is said to be {\it globally
$F$-regular} if there exists an ample line bundle $\mathcal L$
on $X$ such that the section ring 
$$  S(\L) = \oplus_{n \in \mathbb N} {\rm H}^0(X,\mathcal 
L^n),$$
is strongly $F$-regular. 
\end {Definition}  

In case $X$ is globally $F$-regular it can be shown 
that the section ring of any ample line bundle is 
strongly $F$-regular. Moreover, all the local rings 
$\mathcal O_{X,x}$, at points $x$ of $X$, can be 
shown to be strongly $F$-regular (see \cite{SmithFreg}).

Let $X$ be a variety over $k$. The {\it 
absolute Frobenius morphism} on $X$ is the morphism 
$F : X \rightarrow X$ of schemes, which is the identity 
on the set of points and where the associated map 
of sheaves 
$$ F^\sharp : \mathcal O_X  \rightarrow F_* \mathcal O_X$$
is the $p$-th power map. 

Let $D$ be an effective Cartier divisor on $X$ and let $s$ 
denote the associated section of the line bundle $\mathcal O_X(D)$. We 
define $X$ to be {\it Frobenius split along $D$} if the map 
of $\mathcal O_X$-modules
$$ \mathcal O_X \rightarrow F_* \mathcal O_X(D), $$
$$ 1 \rightarrow s,$$
splits. When $D = 0$ we also say that $X$ is {\it Frobenius split} 
\cite{MehtaRamanathan}.
Moreover, if there exists a positive  integer $e$ such 
that  the map of $\mathcal O_X$-modules
$$ \mathcal O_X \rightarrow F_*^e \mathcal O_X(D), $$
$$ 1 \rightarrow s,$$
splits, then we say that $X$ is  {\it stably Frobenius split
along $D$} \cite{SmithFreg}. 

\begin{Proposition}
\label{powers}
Let $X$ denote a variety and let $D$ and $D'$ be 
effective divisors on $X$. If $X$ is stably Frobenius split 
along $D$ then $X$ is stably Frobenius split along $p D$. 
Moreover, if $D' \leq D$ and $X$ is stably Frobenius split 
along $D$ then $X$ is also stably Frobenius split along $D'$. 
\end{Proposition}
\begin{proof}
See Section 3 in \cite{SmithFreg}. 
\end{proof}

Consider now the situation when $X$ is Frobenius split 
(along the zero divisor), and let $Y$ denote a closed 
subscheme of $X$ with sheaf of ideals $\I_Y$.
Let $\phi : F_* \O_X \rightarrow \O_X$
denote the Frobenius splitting. We then say that $Y$ is 
{\it compatibly Frobenius split} if $\phi(F_* \I_Y)
\subseteq \I_Y$. When $Y$ is an effective Cartier 
divisor on $X$ this concept relates to Frobenius splitting 
along $Y$ in the following way. 

\begin{Lemma}
\label{CompSplitting}
Let $Y$ denote an effective Cartier divisor on $X$. If
$Y$ is compatibly Frobenius split in $X$ then $X$ is 
Frobenius split along $(p-1)Y$.
\end{Lemma}
\begin{proof}
Assume that $X$ is Frobenius split and let $\phi : 
F_* \O_X \rightarrow \O_X$ denote 
the associated splitting. The ideal sheaf of $Y$
is isomorphic to $\O_X(-Y)$ and the 
inclusion $ \O_X(-Y) \simeq   \I_Y \subseteq 
\O_X$ is given by multiplication with 
the section $s \in  \O_X(Y)$ defining $Y$. If 
$Y$ is compatibly split then the $p$-power map
$$  \O_X(-Y) \simeq  \I_Y \rightarrow 
F_* \I_Y \simeq F_* \O_X(-Y)$$
splits. Tensoring this map with $\O_X(Y)$ induces
a split map
$$ \O_X \rightarrow (F_*\O_X(-Y))
\otimes \O_X(Y) \simeq F_* \O_X((p-1)Y),$$  
$$ 1 \mapsto s^{p-1},$$
where the isomorphism on the left follows from the 
projection formula and the fact that $F^* \L \simeq \L^{\otimes p}$ 
for every line bundle $\L$ on $X$. This completes 
the proof.
\end{proof}

For a projective variety $X$, global 
$F$-regularity is equivalent to  
$X$ being stably Frobenius split along every 
effective Cartier divisor (cf.~Theorem 3.10(c) in 
\cite{SmithFreg}). This leads to the following result.   

\begin{Lemma}
\label{push down} 
Let $\pi : X \rightarrow Y$ be a morphism of projective 
varieties over $k$ satisfying $\pi_* \mathcal O_X = 
\mathcal O_Y$. If $X$ is globally $F$-regular then 
$Y$ is also globally $F$-regular. In particular, 
if $\pi$ is birational, $X$ is globally F-regular
and $Y$ is normal then $Y$ is globally F-regular.
\end{Lemma}
\begin{proof}
Let $D$ denote an effective Cartier divisor on $Y$, and 
let $s$ denote the corresponding section of the associated
line bundle $\mathcal O_Y(D)$. The pull back of $D$ to $X$ 
will be denoted by $D'$, and the associated section of 
$\mathcal O_X(D')$ is denoted by $s'$. Hence, assuming 
that $X$ is globally
F-regular it follows, by Theorem 3.10 in \cite{SmithFreg}, that there 
exists an integer $e$ such that the morphism of $\mathcal O_X$-modules 
$$ \mathcal O_X \rightarrow F_*^e \mathcal O_X(D'),$$
$$ 1 \mapsto s',$$
splits. Applying the functor $\pi_*$ to this split morphism 
and the assumption $\pi_* \mathcal O_X = \mathcal O_Y$ 
we conclude that the morphism 
$$ \mathcal O_Y \rightarrow F_*^e \mathcal O_Y(D),$$
$$ 1 \mapsto s,$$
of $\mathcal O_Y$-modules splits. The globally $F$-regularity
of $Y$ now follows from Theorem 3.10 in \cite{SmithFreg}. The 
``in particular'' statement follows as $\pi_* \mathcal O_X = 
\mathcal O_Y$ if $\pi$ is birational and $Y$ is normal.
\end{proof}

\section{Schubert varieties are globally $F$-regular}

Let $G$ be a connected and simply connected semisimple linear 
algebraic group $G$ over $k$.
Let $B\subset G$ denote a fixed Borel 
subgroup, $T$ a maximal torus in $B$ and $P$ a parabolic
subgroup containing $B$.   
The $B$-action on the flag variety $G/P$ has finitely many 
orbits $C(w)$ parametrized by left cosets $w\in W/W_P$ in the Weyl 
group $W$ of $G$ with respect to the Weyl group $W_P$ of $P$. We let 
$X(w)$ denote the closure of $C(w)$ in $G/P$. This is the 
Schubert variety corresponding
to $w$. The Weyl group comes with a natural partial order
(the Bruhat order) given by $v\leq w$ if and only if
$X(v)\subseteq X(w)$ for $v, w\in W$. We let $\ell(w)$ denote the
length of the Weyl group element $w\in W$.
For details on the theory of linear algebraic groups we refer to  
\cite{Springer}.

\subsection{Bott-Samelson  varieties}

If $s\in W$ is a simple reflection, the Schubert variety $X(s)$ 
coincides with  the variety $P_s/B\simeq \P^1\subset G/B$ for the minimal 
parabolic subgroup $P_s = B \cup B s B$. 

Let $w=(s_1,s_2,\dots, s_l)$ denote a collection of simple 
reflections in $W$ and let $P_i$ denote the minimal 
parabolic subgroup
associated with $s_i$. The product 
$$P_w=P_1 \times P_2 \times \cdots \times P_l$$ 
comes with a right action of $B^l$ defined as 
$$(p_1, p_2, \cdots, p_l)(b_1,b_2,\cdots,b_l) 
= (p_1 b_1, b_1^{-1} p_2 b_2,\cdots,  b_ {l-1}^{-1} p_l b_l).$$

The quotient $Z(w)=P_w/B^l$ is a smooth projective 
variety of dimension $l$ called a 
Bott-Samelson  variety. 

Fix an integer $1 \leq i \leq l$ and consider the 
closed set of points $Z_i$ in $Z((s_1,\dots,s_l))$
which may be represented by a point in $P_1 \times P_2 
\times \cdots \times P_l$ of the form $(p_1, \dots, 
p_l)$ with $p_i$ equal to the identity element in $G$. 
Then $Z_i$ is a irreducible closed subvariety of 
$Z(s_1,\dots,s_l)$ isomorphic to the 
Bott-Samelson  variety
$Z(s_1,\dots,\hat{s_i}, \dots, s_l)$.
For further details and references on Bott-Samelson varieties 
we refer to \cite{LauTho}.

\begin{Proposition}
\label{F-reg of BSDH}
The Bott-Samelson  variety 
$Z=Z(w)$ is globally $F$-regular.  
\end{Proposition}
\begin{proof} 
As $Z$ is smooth it is enough to prove that $Z$
is stably Frobenius split along an ample divisor
$D$ (see Theorem 3.10 \cite{SmithFreg}). By Theorem 1 
in \cite{MehtaRamanathan} there exists
 a Frobenius splitting of $Z$ compatibly 
splitting the effective divisor $\sum_i Z_i$.
(strictly speaking this result in 
\cite{MehtaRamanathan} only deals with the 
case where $s_1 \cdots s_l$ is a reduced expression.
However, the general case follows in exactly the 
same manner). 

By Lemma \ref{CompSplitting} and Proposition 
\ref{powers} this implies that $X$ is stably
Frobenius split along any divisor of the 
form $\sum_i m_i Z_i$, $m_i \in \mathbb{N}$.
But by Lemma 6.1. in \cite{LauTho}
there exists integers $m_i > 0$ such that the 
divisor $\sum_i m_i Z_i$ is ample. This completes
the proof.
\end{proof}

\subsection{Global $F$-regularity of Schubert varieties}

Consider a Schubert variety $X(w)$ in $G/B$ corresponding
to an element $w$ in the Weyl group $W$. Write $w$ as a
product of simple reflections
$$ w = s_1 s_2 \cdots s_l,$$
with $l$ minimal (i.e. $l = \ell(w)$ the length of $w$).
By the Bruhat decomposition the morphism
$$ P_1 \times P_2 \times \cdots \times P_l \rightarrow X(w),$$ 
$$ (p_1, p_2, \cdots, p_l) \mapsto p_1 p_2 \cdots p_l B.$$
induces a birational morphism
$$Z((s_1,\dots,s_l)) \rightarrow X(w).$$
This leads to the following result.

\begin{Theorem}\label{CorollaryFreg}
Let $P$ denote a parabolic subgroup in $G$ which contains
$B$. A Schubert variety $X$ in $G/P$ is 
globally F-regular. 
\end{Theorem}
\begin{proof}
It is well known that Schubert varieties 
are normal. The case $P=B$ follows immediately 
from Lemma \ref{push down} and Proposition 
\ref{F-reg of BSDH}, together with the existence 
of the birational map 
$$Z(s_1,\dots,s_l) \rightarrow X(w),$$
described above.

Consider now a general parabolic subgroup $P$.
The inverse image of $X$, by the canonical map
$\pi : G/B \rightarrow G/P$, is then a Schubert
variety $X(w)$ in $G/B$ for some $w\in W$. 
Choose $w'$ to be a minimal length representative 
for the left coset $w W_P$ in $W/W_P$. By the 
Bruhat decomposition it follows that 
the induced map :
$$ \pi : X(w') \rightarrow X,$$
is birational. As $X(w')$ is globally F-regular
the statement is now a consequence of Lemma 1.2.
\end{proof}

\section{$F$-rationality and $\D$-modules}

Let $R$ denote a commutative algebra over a perfect field $k$.
The ring of $k$-linear differential operator $D_k(R)$ on 
$R$ is an $R \otimes_k R$-subalgebra of ${\rm End}_k(R)$
defined by 
$$  D_k(R) = \{ \phi \in {\rm End}_k(R) : 
I^ n \cdot \phi = 0, ~n \gg 0 \}, $$ 
where $I$ denotes the kernel of the product map 
$R \otimes_k R \rightarrow R$. The $R\otimes_k R$-submodules 
$$ D_k^n (R) = \{ \phi \in {\rm End}_k(R) : 
I^ {n+1} \cdot \phi = 0\}, $$ 
defines a filtration of $D_k(R)$. Elements in 
$D_k^n (R)$ is called differential operators of
degree $\leq n$. When $I$ is a finitely generated 
ideal there is a second filtration of $D_k(R)$  
given by the $R\otimes_k R$-submodules 
$$ D_{k}^{(n)}(R) = \{ \phi \in {\rm End}_k(R) : 
I^{(n+1)} \cdot \phi = 0\}, $$ 
where $I^{(n)}$ denotes the ideal in $R\otimes_k R$
generated by elements of the form $a^n$, $a \in I$. 
This filtration is particularly nice when the characteristic
$p$ of $k$ is positive. In this case $I^{(p^n)}$ is 
generated  by elements of the form 
$ a^{p^n} \otimes 1 - 1 \otimes a^{p^n}, $
and hence $D_k^{(p^n-1)} (R) = {\rm End}_{R^{p^n}} (R),$
where $R^{p^n}$ denotes the subring of $R$ of $p^n$-powers
(here we use that $k$ is algebraically 
closed and hence perfect).
In particular, 
$$ D_k(R) = \bigcup_n {\rm End}_{R^{p^n}} (R).$$
The right side of this equation shows that $D_k(R)$ 
is independent of $k$, and we therefore suppress  
$k$ from the notation and write $D(R)$ instead 
of $D_k(R)$.

\begin{Lemma}
\label{ring}
Assume that $k$ has positive characteristic $p$
and that $R$ is a finitely generated $k$-algebra.
For every multiplicative subset $S$ of $R$ there
exists a natural isomorphism of left $R_S$-modules
$$(D(R))_{S} \simeq D(R_{S}),
$$
where the localization on the left is performed
as a left $R$-module.
\end{Lemma}
\begin{proof}
Fix a positive integer $n$. As $R$ is a finitely generated
$k$-algebra it is finitely generated as a module over the 
subring $R^{p^n}$. This implies that the exists a natural
isomorphism 
$${\rm End}_{R^{p^n}} (R)_S \simeq {\rm End}_{R_S^{p^n}} (R_S),$$
Now conclude the argument by using the description of $D(R)$ 
above (in positive characteristic).
\end{proof}

\subsection{Sheaves of differential operators}

Let $X$ be a variety over $k$. The sheaf of
$k$-linear differential operators $\D_X$ 
on $X$ (cf.~\cite{EGA}, \S 16.8) is a 
$\O_X$-bisubalgebra of ${\mathcal End}_k(\O_X)$ 
which is quasicoherent for both $\O_X$-modules 
structures. When $X ={\rm Spec}(R)$ is affine,
the sheaf 
$\D_X$ coincides with the quasicoherent 
$\O_X$-bialgebra associated to the $R$-bialgebra 
$D_k(R)$ defined above.


\begin{Lemma}\label{LemmaLocD}
Assume that $k$ has positive characteristic.
For every $x\in X$ the natural morphism of algebras
$$
\phi : (\D_X)_x  \rightarrow D(\O_{X,x})
$$
is an isomorphism as $\O_{X, x}$-bimodules.
\end{Lemma}
\begin{proof}
We may assume that $X = {\rm Spec}(R)$ is affine.
When  $\mathfrak p$ denotes the prime ideal in $R$
associated with $x$ we may identify $\O_{X,x}$ with
$R_\mathfrak p$ and $(\D_X)_x$ with $D(R)_\mathfrak p$ 
(as a left $\O_{X,x}$-bimodule). The result now follows
from Lemma \ref{ring} above.
\end{proof}

\subsection{Local cohomology and $\D$-modules}
A sheaf of abelian groups $\F$ on a variety $X$ is 
called a $\D_X$-module, if $\F$ is a module over the
sheaf of algebras $\D_X$ such that the $\O_X$-structure, 
induced by the inclusion of $\O_X$ in $\D_X$, 
is quasicoherent.

For a locally closed subset $C$ of  $X$ 
and a sheaf $\F$ of abelian groups on $X$, we let 
$\SH^i_C(\F)$ denote the $i$-th local cohomology sheaf 
with support in $C$ (cf.~\cite{LC}).  If $\F$ is a 
$\D_X$-module, then $\SH^i_C(\F)$ is a $\D_X$ module 
for any locally closed subset $C\subseteq X$.
In particular, we may regard $\SH^i_C(\O_X)$
as a $\D_X$-module. 

\begin{Lemma}
\label{finite length}
Let  $C$ be a locally closed subset of a smooth variety 
$X$. Then the $\D_X$-module $\SH^i_C(\O_X)$ has finite 
length.
\end{Lemma}
\begin{proof}
The proof depends on the concept of filtration 
holonomicity defined by B\" ogvad in  \cite{Bogvad}.
By Prop.3.7. in \cite{Bogvad} the $\D_X$-module 
$\SH^i_C(\O_X)$ is filtration holonomic. But 
any filtration holonomic module has finite 
length  as a 
$\D_X$-module (Prop.3.2. in  \cite{Bogvad}).
\end{proof}

\subsection{Support of finite length $\D$-modules}

Let $X$ be a variety over $k$. In this section we prove that 
the support of any $\D_X$-module of finite length is closed. 

\begin{Lemma}
\label{restriction is simple}
Let $X$ be a variety over $k$ and let $\F$ denote a simple
$\D_X$-module. Then $\F_{|U} $ is a simple $\D_U$-module 
for any open subset  $U$ of $X$.
\end{Lemma}
\begin{proof}
Let $i : U \rightarrow X$ denote the inclusion map. The  
restriction map $\F \rightarrow i_* \F_{|U}$ is a map of 
$\D_X$-modules. Let $\M$ denote a $\D_U$-submodule 
of $\F_{|U}$ and consider the composed map
$$ \phi : \F \rightarrow i_* \F_{|U} \rightarrow i_* ( \F_{|U}/ \M).$$
The kernel ker$(\phi)$ is a $\D_X$-submodule of $\F$ and hence 
either ker$(\phi)$ equals  $\F$ or $0$. In particular,
the restriction of ker$(\phi)$ to $U$ is either $\F_{|U} $
or $0$. But, ker$(\phi)|_U$
equals $\M$. This completes the proof.
\end{proof}

\begin{Lemma}
\label{closed support}
Let $X$ be a variety over $k$ and let $\F$ denote a 
$\D_X$-module of finite length. Then the support of 
$\F$ is closed.
\end{Lemma}
\begin{proof}
As $\F$ has finite length there exists a filtration 
$$ 0 = \F_0 \subseteq \F_1 \subseteq 
\cdots \subseteq \F_m = \F,$$
by $\D_X$-submodules, such that the quotients $\L_i = \F_i/
\F_{i-1}$ are simple $\D_X$-modules. Moreover, the support
of $\F$ is the union of the supports of $\L_i, i=1, \dots 
m$. This reduces the statement to the case when $\F$ is 
a simple $\D_X$-module. So assume now that $\F$ is simple.

By Lemma \ref{restriction is simple} we may furthermore 
assume that $X$ is affine. Hence, there exists a global 
(nonzero) section $s$ of $\F$. By simplicity we must have 
$\D_X \cdot s = \F$, and hence the support of $\F$ 
coincides with the support of $s$ which is closed.
\end{proof}

\subsection{$F$-rationality and $\D$-modules}

The concept of $F$-rationality comes from the theory of
tight closure in commutative algebra. A local commutative
ring in positive characteristic is called $F$-rational
if every parameter ideal in $R$ is tightly closed. For the
definition of tight closure see \cite{HH}. One may prove
that every ideal in a strongly $F$-regular ring is tightly
closed. In particular it follows that strongly $F$-regular
rings are $F$-rational.

We quote the following crucial result by M.~Blickle.

\begin{Theorem}[Corollary 4.10 in \cite{Blickle}] \label{TheoremBlickle}
Let $R$ be regular, local and $F$-finite.
Let $A = R/I$ be a domain of codimension $c$. If $A$ is 
$F$-rational, then $H^c_I(R)$ is $D_R$-simple.
\end{Theorem} 

Here $F$-finite means that $R$ is finitely generated as a 
module over the subring $R^p$ of $p$-th powers. In particular,
the Theorem applies when $R$ is the localization of a finitely
generated $k$-algebra. This leads to the following global
result.

\begin{Proposition}
\label{simplicity}
Let $Y$ be an irreducible closed subvariety of codimension $c$
of a smooth variety $X$ over $k$. If all the local rings $\O_{Y,y}$,
$y \in Y$, are $F$-rational, then the $\D_X$-module $\SH_Y^c(\O_X)$ 
is simple.    
\end{Proposition} 
\begin{proof}
By Lemma \ref{finite length} there exists a finite length submodule 
$\M$ of $\SH_Y^c(\O_X)$ such that the associated 
quotient $\L$ is a (nonzero) simple  $\D_X$-module. Consider 
the corresponding short exact sequence: 
$$ 0 \rightarrow \M \rightarrow \SH_Y^c(\O_X) \rightarrow 
\L \rightarrow 0.$$
The induced map of stalks 
$$ 0 \rightarrow \M_x \rightarrow \SH_Y^c(\O_X)_x \rightarrow 
\L_x \rightarrow 0.$$
at a point $x$ in $X$, is  then a short exact sequence 
of $\D_{X,x}$-modules. Using Theorem \ref{TheoremBlickle} and 
Lemma \ref{LemmaLocD} we conclude that the middle term is a
simple $\D_{X,x}$-module and hence that either $\L_x$ or $\M_x$ 
is zero. In particular, the support Supp$(\SH_Y^c(\O_X))$,
which is $Y$,  is the disjoint union of the 
supports of $\M$ and $\L$. But, by Lemma \ref{closed support},
the supports of $\L$ and $\M$ are closed, and as $Y$ is irreducible
we conclude that either Supp$(\M)$ or Supp$(\L)$ is empty.
Hence, $\M=0$ (as $\L \neq 0$ by 
assumption) and $\SH_Y^c(\O_X)$ thereby coincides with 
the simple module $\L$.
\end{proof}

\section{$\D$-modules on $G/B$}

Let $X$ denote a (generalized) flag variety $G/P$ and 
let $Y$ denote a Schubert variety in $X$ of codimension 
$c$. We assume that 
the characteristic of the ground field $k$ is positive.

\begin{Theorem}
\label{simple}
The $\D_X$-module $\SH^c_{Y}(\O_X)$ is simple.
\end{Theorem}
\begin{proof}
By Proposition \ref{simplicity} it is enough to prove 
that the local rings $\O_{Y,y}, y \in Y,$ are $F$-rational.
But $Y$ is globally $F$-regular by Theorem \ref{CorollaryFreg}. 
In particular, the local rings $\O_{Y,y}$ are strongly 
$F$-regular and hence $F$-rational.
\end{proof}

For a smooth algebraic variety $X$ with an action
of an algebraic group $H$, there is a natural notion of 
an $H$-equivariant $\D_X$-module (see \cite{KAS}, p.~82).

Consider the category of $B$-equivariant holonomic
$\D_X$-modules for $X=G/B$. Independently of $\kar(k)$
a simple module in this category is uniquely given by
its support which is a Schubert variety (see \cite{Bogvad}, 
Theorem 4.6 for the positive characteristic case and \cite{BRK}, 
\S6.4 for the characteristic zero case). We let $\L(w)$ denote 
the simple module with support $X(w)$ for $w\in W$. Independently
of $\kar(k)$ one may prove that $\L(w) \subseteq \SH^c_{X(w)}(\O_X)$ 
with equality if $X(w)$ is
smooth.

As a consequence 
of Theorem \ref{simple} we have the following result.

\begin{Theorem}\label{TheoremXw}
Suppose that $X$ is over a field of positive characteristic.
The unqiue simple $B$-equivariant and holonomic $\D_X$-module 
$\L(w)$ with support $X(w)$ is 
isomorphic to the local cohomology module
$$
\SH^c_{X(w)}(\O_X)
$$
where $c$ denotes the codimension of $X(w)$.
\end{Theorem}

In the next section we give an example showing that Theorem
\ref{TheoremXw} does not hold in characteristic zero.

\subsection{Non-simplicity in characteristic zero}

Consider the algebraic group $G=\SL_4$ over the field $k$ 
with $\kar (k)=0$ along with its $6$-dimensional 
flag variety $X = G/B$. The Weyl group $W$ of $G$ is 
generated by the
simple reflections $s_1,s_2,s_3$ numbered from left to right
in the Dynkin diagram. We let $P_{v,w}$ 
denote the Kazhdan-Lusztig polynomial associated with $v, w\in W$ 
(cf.~\cite{KLU}). 

Consider $w=s_1 s_2 s_3 s_2 s_1$. Then $X(w)$ is a
codimension one Schubert variety in $X$. One may check 
that (in fact the singular locus of $X(w)$ is
$X(s_1 s_3)$) 
$$
P_{v, w}(q) =
\begin{cases}
1+q & \text{if $v\leq s_1s_3$}\\
1 & \text{if $v\nleq s_1 s_3$}.
\end{cases}
$$ 
According to the Kazhdan-Lusztig conjecture (proved in \cite{BB},\cite{BRK},\cite{KLU}) we have
the following formula for the simple module $\L(w)$
in the Grothendieck group of the category of $B$-equivariant and
holonomic $\D$-modules on $X$
$$
[\L(w) ]=\sum_{v\leq w}(-1)^{\ell(w)-\ell(v)}P_{v,w}(1)[\SH
^{6-l(v)}_{C(v)}(\O_X )]. 
$$
The Schubert variety 
$X(w)$  is a local complete
intersection since it has codimension one. Therefore   
$\SH^j_{X(w)}(\O_X)$ is non-vanishing if and only if $j=1$. 
It follows by~\cite{KL} (see \S \ref{GC}) that we have 
an exact sequence 
$$
0\rightarrow \SH^1_{X(w)}(\O_X)\rightarrow
\SH^1_{C(w)}(\O_X)\rightarrow \bigoplus _{v\leq w,
  \ell(w)-\ell(v)=1}\SH^2_{C(v)}(\O_X)\rightarrow \dots.    
$$
From this we get the formula 
\begin{eqnarray*}
& & [\SH^1_{X(w)}(\O_X)]=\sum _{v\leq w}(-1)^{\ell(v)-\ell(w)}[\SH
  ^{6-\ell(v)}_{C(v)}(\O_X )]= \\ 
& &  [\mathcal{L}(w)]+\sum _{v< w}(-1)^{\ell(v)-5}(1-P_{v,w}(1))[\SH
 ^{6-\ell(v)}_{C(v)}(\O_X )]= \\
& &  [\mathcal{L}(w)] +\sum _{v\leq
  s_1 s_3}(-1)^{\ell(v)-2}[\SH ^{6-\ell(v)}_{C(v)}(\O_X)]= \\
& & [\mathcal{L}(w)] +[\mathcal{L}(s_1 s_3)] \\
\end{eqnarray*}
in the Grothendieck group. It follows that
$\SH^1_{X(w)}(\O_X)$ is not a simple $\D_X$-module.
\label{nonsmpl}

\subsection{Decomposition of dual Verma modules}
\label{GC}
In this section we consider the category of $B$-equivariant and
(filtration) holonomic $\D_X$-modules on $X = G/B$ over a field
of positive characteristic.
Associated to the $B$-orbit $C(w)\subset G/B$ we have the
local cohomology sheaf $\M(w) = \SH^c_{C(w)}(\O_X)$, where $c$ denotes
the codimension of $C(w)$. The simple module 
$\L(w) = \SH^c_{X(w)}(\O_X)$ admits a natural injection
into $\M(w)$.
Over a field 
of characteristic zero the decomposition of $\M(w)$ into simple modules 
ultimately rests on deep
properties of intersection cohomology. The situation in positive
characteristic is quite different. 
Applying the functor $\underline{\Gamma}_{X(w)}(-)$ 
to the local Grothendieck-Cousin 
complex for $\O_X$ we get a complex (see \cite{KL})
\begin{equation}\label{exact}
0\rightarrow \L(w) \rightarrow \M^0 \rightarrow \M^1 \rightarrow
\cdots \rightarrow \M^i \rightarrow \cdots,
\end{equation}
where 
$$
\M^i = \bigoplus_{y\leq w, \ell(y) = \ell(w) - i} \M(y)
$$
and $\L(w)$ sits in degree $c-1$.
The complex $\M^\bullet$ in (\ref{exact}) computes the cohomology
$\SH^j_{X(w)}(\O_X)$. Since $X(w)$ is Cohen-Macaulay it follows by
(\cite{PS}, Chapitre III, Proposition (4.1)) that 
(\ref{exact}) is exact. This shows that we have the formula
\begin{equation}\label{invert}
[\L(w)] = \sum_{y\leq w} (-1)^{\ell(w) - \ell(y)} [\M(y)]
\end{equation}
in the Grothendieck group
of holonomic and $B$-equivariant $\D$-modules on $G/B$ in positive
characteristic. Using Verma's identity (cf.~\cite{KLU}, Remarks 3.3(b))
$$
\sum_{x\leq z \leq y} (-1)^{\ell(x)}(-1)^{\ell(z)} = \delta_{x, y}
$$
for all $x\leq y$ in $W$, one can 
invert (\ref{invert}) to get the decomposition
$$
[\M(w)] = \sum_{y\leq w} [\L(y)]
$$
for the ``dual Verma module'' $\M(w)$. This shows that the 
decomposition of $\M(w)$ into simple modules is multiplicity free and
that $[\M(w):\L(y)] = 1$ if $y\leq w$ and $0$ otherwise.
In the characteristic zero setting this decomposition is 
given by the value of inverse Kazhdan-Lusztig polynomials at $1$. 
The authors are unaware of any natural correspondence between the
categories of holonomic and equivariant $\D$-modules in zero and
positive characteristic even though they share many properties
(e.g.~simple modules and dual Verma modules parametrized by Schubert
varieties). 

\bibliographystyle{amsplain}

\begin{thebibliography}{10}

\bibitem{BB}
A.~Beilinson and J.~Bernstein,
Localisation de $\g$-modules,
\emph{C.~R.~Acad.~Sci. Paris}~{\bf 292}(1981), 15--18.

\bibitem{Blickle}
Manuel Blickle, \emph{The intersection homology {$D$}-module in finite
  characteristic}, Math.~Ann., (to appear).

\bibitem{Bogvad}
R.~B\"ogvad, \emph{Some results on {$D$}-modules on {B}orel varieties},
  J.~Algebra \textbf{173} (1995), 638--667.

\bibitem{Bogvad02}
\bysame, \emph{An analogue of holonomic {$D$}-modules on smooth varieties in
  positive characteristics}, Homology Homotopy Appl. \textbf{4} (2002),
  83--116.

\bibitem{BRK}
J.~L. Brylinski and M.~Kashiwara, \emph{Kazhdan-Lusztig conjecture and
  holonomic systems}, Invent.~Math. \textbf{64} (1981), 387--410.

\bibitem{EGA}
A.~Grothendieck, \emph{{EGA IV, Etude locale des schemas et des morphismes de
  schemas}}, Inst.~Hautes \'Etudes Sci.~Publ.~Math. \textbf{32} (1967)

\bibitem{LC}
\bysame, \emph{Local {C}ohomology}, Springer-Verlag, 1966.


\bibitem{EK}
M.~Emerton and M.~Kisin, \emph{A Riemann-Hilbert correspondence for
  unit $F$-crystals}, Ast\'erisque \textbf{293} (2004). 


\bibitem{HH}
M.~Hochster and C.~Huneke, \emph{{Tight closure and strong $F$-regularity}},
  M\'emoires de la Soci\'et\'e Math\'ematique de France \textbf{38} (1989),
  119--133.

\bibitem{KAS}
M.~Kashiwara, \emph{Representation theory and {$D$}-modules on flag varieties},
  Ast\'erisque \textbf{173--174} (1989), 55--109.

\bibitem{KL}
M.~Kashiwara and N.~Lauritzen, \emph{Local cohomology and {$\D$}-affinity in
  positive characteristic}, C.~R.~Acad.~Sci.~Paris \textbf{335} (2002),
  993--996.

\bibitem{KLU}
D.~Kazhdan and G.~Lusztig, \emph{{Representations of Coxeter groups and Hecke
  algebras}}, Invent.~Math. \textbf{53} (1979), 165--184.

\bibitem{LM}
V.~Lakshmibai and P.~Magyar, \emph{Degeneracy schemes and Schubert varieties}, 
Int.~Math.~Res.~Notices \textbf{12} (1998), 627--640.

\bibitem{LauTho}
N.~Lauritzen and J.~F.~Thomsen, \emph{{Line bundles on Bott-Samelson 
varieties}}, J.~Alg.~Geom. \textbf{13} (2004), 461--473.

\bibitem{MehtaRamanathan}
V.~Mehta and A.~Ramanathan \emph{Frobenius splitting and cohomology 
vanishing for Schubert varieties},
  Ann.~Math. \textbf{122} (1985), 27--40.

\bibitem{Ram}
A.~Ramanathan, \emph{Equations defining Schubert varieties and Frobenius splitting
of diagonals}, Pub.~ Math.~ IHES  \textbf{65} (1987), 61--90.

\bibitem{PS}
C.~Peskine, L.~Szpiro, \emph{Dimension projective finie et cohomologie locale}, Pub.~ Math.~ IHES  \textbf{42} (1973), 47--119.

\bibitem{SmithFreg}
K.~E.~Smith, \emph{{Globally $F$-regular varieties: Applications to
  vanishing theorems for quotients of Fano varieties}}, Mich.~J.~Math.
  \textbf{48} (2000), 553--572.

\bibitem{Springer}
T.~Springer, \emph{Linear Algebraic Groups}, Springer Verlag (2nd
edition), 1998

\end{thebibliography}

\providecommand{\bysame}{\leavevmode\hbox to3em{\hrulefill}\thinspace}
\providecommand{\MR}{\relax\ifhmode\unskip\space\fi MR }
\providecommand{\MRhref}[2]{%
  \href{http://www.ams.org/mathscinet-getitem?mr=#1}{#2}
}
\providecommand{\href}[2]{#2}

\end{document}